\documentclass[11pt]{article}
\usepackage[english]{babel}
\usepackage{amsfonts, amsmath, amsthm, amssymb,amscd,indentfirst}
\usepackage{graphicx}
\usepackage{wrapfig}
\usepackage[rflt]{floatflt}
\usepackage{amsmath,amssymb,latexsym,indentfirst}
\usepackage{enumerate,pst-all}
\usepackage{makeidx}
\usepackage{color}
\usepackage{times}
\usepackage{palatino}
\newtheorem{theorem}{Theorem}[section]

\newtheorem{lemma}{Lemma}[section]

\newtheorem{remark}{Remark}[section]

\def\Cc{{\cal C}}

\def\Cc{{\cal C}}
\def\R{{\cal R}}

\def\Uc{{\cal U}}

\def\Sc{{\cal S}}

\def\Lc{{\cal L}}
\def\Gc{{\cal G}}
\def\Fc{{\cal F}}
\def\Rc{{\cal R}}

\def\Ac{{\cal A}}

\def\Cc{{\cal C}}
\def\gm{{g_{\mu}}}

\def\Mc{{\cal M}}

\def\Lc{{\cal L}}
\def\Gc{{\cal G}}
\def\Fc{{\cal F}}
\def\Rc{{\cal R}}

\def\Ac{{\cal A}}

\def\Cc{\mathcal{C}}

\def\Rc{{\cal R}}
\def\Lc{{\cal L}}
\def\Fc{{\cal F}}

\def\Uc{{\cal U}}

\def\Mc{{\cal M}}

\begin{document}

\title{The Yamabe problem for Gauss-Bonnet curvatures: a local result around space forms}

\author{Levi Lopes de Lima\thanks{Federal University of Cear\'a,
Department of Mathematics, Campus do Pici, R. Humberto Monte, s/n, 60455-760,
Fortaleza/CE, Brazil ({\tt levi@mat.ufc.br}). Partially supported by CNPq.}
 \and Newton Lu\'is Santos\thanks{Federal University of
Piau\'{\i},
 Department of Mathematics, Campus
Petronio Portela, Ininga, 64049-550 Teresina/PI, Brazil ({\tt newtonls@ufpi.br}).
Partially supported by a CNPq Posdoctoral Grant.}
} \maketitle

\begin{abstract}
It is shown in \cite{L2} that metrics with {\em constant} $2k$-Gauss-Bonnet curvature on a closed $n$-dimensional manifold, $2\leq 2k <n$, are critical points for a certain Hilbert type functional  with respect to volume preserving conformal variations. This motivates the corresponding Yamabe problem: is it true that any metric on a closed manifold is conformal to a metric with constant $2k$-Gauss-Bonnet curvature? Using perturbative methods we affirmatively answer this question for small perturbations of certain space forms. More precisely, if $(X,g)$ is a
non-flat closed space form
{\em not} isometric to a round sphere
we show the existence of a neighborhood $U$ of $g$ in the space of metrics such that any $g'\in U$ is conformal to a metric whose $2k$-Gauss-Bonnet curvature is constant. \end{abstract}

\section{Introduction and statement of results}

Let $X$ be a connected manifold of dimension $n\geq 3$. We  denote by $\Omega^p(X)$ the
space of exterior $p$-forms on $X$. The space of double forms
of bi-degree $(p,q)$ is
$$
\Omega^{p,q}(X)=\Omega^p(X)\otimes_{\Omega^0(X)}\Omega^q(X),
$$
so that the direct sum
$\Omega(X)=\oplus_{p,q\geq 0}\Omega^{p,q}(X)$ is a bi-graded associative algebra.
For example, any bilinear form on tangent vectors is a
$(1,1)$-form. In particular, if $\Mc(X)$ is the space of Riemannian metrics on $X$ then any $g\in\Mc(X)$ is a
$(1,1)$-form and the curvature tensor $R_g$ of $g$ is a
$(2,2)$-form. In fact, if we define $\Cc^p(X)\subset \Omega^{p,p}(X)$ to be the space of $(p,p)$-forms satisfying the symmetry condition
$$
\omega(v_1\wedge\ldots\wedge v_{p}\otimes w_1\wedge\ldots\wedge
w_{p})=\omega(w_1\wedge\ldots\wedge w_{p}\otimes v_1\wedge\ldots\wedge
v_{p}),
$$
then any bilinear symmetric form on tangents vectors ($g$ in particular) belongs   to $\Cc^1(X)$, and moreover $R_g\in\Cc^2(X)$. An account of the formalism of double forms used in this paper can be found in \cite{L1}, \cite{L2}.

Let us consider a Riemannian manifold $(X,g)$ with $X$ as above. Then multiplication
by the metric defines a map $g:\Omega^{p-1,q-1}_x(X)\to\Omega^{p,q}_x(X)$ whose adjoint with respect to the natural inner product
$\langle\,,\rangle$ on $\Omega_x^{*,*}(X)$ induced by $g$ is the
contraction operator $c_g:\Omega^{p,q}_x(X)\to\Omega^{p-1,q-1}_x(X)$ defined by
$$
c_g\omega(v_1\wedge\ldots\wedge v_{p-1}\otimes w_1\wedge\ldots\wedge
w_{p-1})=\sum_i \omega(e_i\wedge v_1\wedge\ldots\wedge
v_{p-1}\otimes e_i\wedge w_1\wedge\ldots\wedge w_{p-1}),
$$
where $\{e_i\}$ is a local orthonormal tangent frame.

With this formalism at hand, we set, for any $1\leq k\leq n/2$,
\begin{equation}\label{def}
\Sc^{(2k)}_g=\frac{1}{(2k)!}c^{2k}_gR^k_g.
\end{equation}
This is the $2k$-{\em Gauss-Bonnet curvature} of $g$. Thus, up to a constant this invariant interpolates between the scalar curvature ($k=1$) and  the Gauss-Bonnet integrand ($2k=n$).
Locally,
\begin{equation}\label{fully}
\Sc^{(2k)}_g=c_{n,k}\sum\delta_{i_1\ldots i_{2k}}^{j_1\ldots j_{2k}}R_{j_1j_2}^{i_1i_2}\ldots R_{j_{2k-1}j_{2k}}^{i_{2k-1}i_{2k}},
\end{equation}
where $R_{ij}^{kl}$ are the components of $R_g$ in an orthonormal frame. Thus, $\Sc^{(2k)}_g$ is  homogeneous of degree $k$ in $R_g$. If $X$ is isometrically embedded in a Euclidean space as a hypersurface then $\Sc^{(2k)}_g$ is, up to a constant, the elementary symmetric function of order $2k$ in the principal curvatures. Also, the metric invariants $\int_X\Sc^{(2k)}_g\nu_g$, where $\nu_g$ is the volume element of $g$, appear both in Weyl's formula for the volume of tubes \cite{G} and  Chern's kinematic formula \cite{Ch}, besides playing a key role in Stringy Gravity \cite{C-B}.

\begin{remark}\label{const}
Notice that in general we have
\begin{equation}\label{diff}
\Sc^{(2k)}_g=\frac{1}{(2k)!}c_g\Rc_g^{(2k)},
\end{equation}
where
\begin{equation}\label{ricci}
\Rc^{(2k)}_g=c^{2k-1}_gR^k_g
\end{equation}
is the $2k$-{\em Ricci tensor} of $g$ (again for $k=1$ we recover the standard notion of Ricci tensor).
Moreover, using the language of double forms, that a Riemannian manifold $(X,g)$ has
constant sectional curvature $\mu$ can be characterized  by the identity
$R_g=\frac{\mu}{2}g^2$, where $g=g_\mu$ is the constant curvature metric. In this case, the $2k$-Ricci tensor and the $2k$-Gauss-Bonnet curvature are, respectively,
\begin{equation}\label{lambda}
\R^{(2k)}_{\gm}=\frac{(n-1)!(2k)!}{(n-2k)!2^k}\mu^k\gm,\qquad \mathcal \Sc^{(2k)}_{\gm}=\frac{n!}{(n-2k)!2^k}\mu^k.
\end{equation}
\end{remark}

Assume from now on that $X$ is {\em closed}. In this setting it is proved in \cite{L2} that for $2\leq 2k<n$ metrics with {\em constant} $2k$-Gauss-Bonnet curvature are critical for the so-called {\em Hilbert-Einstein-Lovelock} functional
\begin{equation}\label{hillove}
g\in\Mc(X)\mapsto \int_X\Sc^{(2k)}_g\nu_g\in\mathbb{R}
\end{equation}
with respect to volume preserving conformal variations; for $k=1$ this is a classical result \cite{B}. As suggested in \cite{L2}, this motivates a sort of Yamabe problem for Gauss-Bonnet curvatures: {\em is it true that any metric on a closed manifold is conformal to a metric with constant $2k$-Gauss-Bonnet curvature?}

For $k=1$ this is a classical problem in Conformal Geometry solved by Schoen after previous contributions by Yamabe, Trudinger and Aubin \cite{LP}. If $k\geq 2$
and $(X,g)$ is locally conformally flat, the Yamabe problem for $\Sc_g^{(2k)}$ is equivalent to the so-called $\sigma_k$-Yamabe problem \cite{V}, which consists of finding conformal metrics whose elementary symmetric function of order $k$ of the Schouten tensor is constant; see \cite{L3}. For this class of manifolds the Yamabe problem has been solved in the affirmative in \cite{GW} and \cite{LL}, assuming that the background metric satisfies a certain ellipticity condition. Our main result, Theorem \ref{main2} below, adds a new class of manifolds to this list, namely, small perturbations of certain space forms, and  seems to provide the first examples of non-locally conformally flat metrics for which the Yamabe problem for the Gauss-Bonnet curvatures is solved.

More precisely, let $\Fc_n$ be the class of $n$-dimensional non-flat closed space forms {\em not} isometric to a round sphere. Moreover, if $(X,\gm)\in\Fc_n$, $\Omega^0_+(X)$ is the set of smooth positive functions on $X$ and $1:X\to\mathbb{R}$ is the function identically equal to $1$. With this notation at hand the following result is obtained.

\begin{theorem}\label{main2}
Assume $4\leq 2k<n$ and let $(X,\gm)\in\Fc_n$ with ${\rm vol}(X,\gm)=\nu$. Then the space $\Mc^{k}_{\nu}(X)$ of metrics on $X$ with constant $2k$-Gauss-Bonnet curvature and volume $\nu$ has, around $\gm$, the structure of an (infinite dimensional) ILH-submanifold of $\Mc(X)$. Moreover, the map $\xi:\Omega^0_+(X)\times\Mc_{\nu}^{k}(X)\to\Mc(X)$, given by $\xi(f,g)=fg$, is ILH-smooth around $(1,\gm)$ and its derivative at $(1,\gm)$ is an isomorphism. In particular, there exists a neighborhood $U$ of $\gm$ in $\Mc(X)$ such that any metric in $U$ is conformal to a metric with constant $2k$-Gauss-Bonnet curvature and volume $\nu$.
\end{theorem}

For the $ILH$-terminology we refer to \cite{O}. We also mention that the proof leading to Theorem \ref{main2} is  modeled upon an argument due to Koiso \cite{K}, where a local Yamabe type result  has been proved  in case $k=1$ for a  class of Riemannian manifolds containing $\Fc_n$.

\begin{remark}\label{fullynon}
The result in Theorem \ref{main2} amounts to solving a second order {\em fully nonlinear} equation for the conformal factor. We refer to \cite{L3} where the equation is displayed for $k=2$.
\end{remark}

\begin{remark}\label{sphere}
Note that the local Yamabe result, as formulated  in Theorem \ref{main2}, does not hold for $(X,\gm)$ a round sphere. In fact, by pulling back the standard metric by the flow of a conformal field we obtain a one parameter family of conformal metrics with the same constant scalar curvature and volume.
\end{remark}

\section{Linearizing the Gauss-Bonnet curvatures}

The key step in the proof of Theorem \ref{main2} is to compute the linearization $\dot\Sc^{(2k)}_g:\Cc^1(X)\to \Omega^0(X)$ of the operator $g\in\Mc(X)\mapsto {\mathcal S}_g^{(2k)}\in\Omega^0(X)$ at a metric $\gm$ with constant sectional curvature $\mu\neq 0$. For a general metric $g$ we get from (\ref{diff})
that
\begin{equation*}\label{lin1}
\dot{\mathcal S}_g^{2k}h=\frac{1}{(2k)!}\left((\dot c_gh){\mathcal R}_g^{(2k)}+c_g
\dot{\mathcal R}_g^{(2k)}h\right).
\end{equation*}
Notice  that the first term in the right-hand side above is of order zero in $h$ (no derivatives) while the second one has derivatives up to second order of $h$ and hence determines the analytic properties of the linearization (its symbol). Since from (\ref{ricci}) we have
$$
\dot\Rc_g^{(2k)}h=A_gh+B_gh,
$$
where
\begin{equation*}\label{lin3}
A_gh=(2k-1)(\dot c_gh)c_g^{2k-2}R_g^k
\end{equation*}
and
\begin{equation*}\label{lin4}
B_gh=kc_g^{2k-1}R^{k-1}_g\dot R_gh,
\end{equation*}
we see that, for $k\geq 2$, this higher order term  depends on the {\em full} curvature tensor and not only on the Ricci tensor as it is the case if $k=1$.
Thus, in order to get a workable expression for $\dot\Sc_g^{(2k)}$, $k\geq 2$, one is naturally led to assume that $R_g$ should depend on $g$ in a simple manner.
If $\gm$ has constant sectional curvature $\mu\neq 0$, so that $R_\gm=\mu g_{\mu}^2/2$ (see Remark \ref{const}), it is shown in \cite{dLS}  that
$$
A_\gm h=
(2k-1)(n-2)C_{n,k}\mu^k(h-({\mathsf {tr}}_\gm h)\gm)
$$
and
\begin{eqnarray*}
B_\gm
h & = &  C_{n,k}\mu^{k-1}\bigg((k-1)\left(\Delta_\gm \mathsf {tr}_\gm h+\delta_\gm (\delta_\gm h)+ (n-1)\mu \mathsf {tr}_{g_\mu} h\right)\gm +  \nonumber\\
& &\qquad +\frac{n-2k}2 \Big(\nabla^*\nabla h-
\nabla d\mathsf {tr}_\gm h-\\
& & \qquad\qquad-2\delta^{*}_\gm\delta_\gm h+2(n-1)\mu h \Big)\bigg),\nonumber
\end{eqnarray*}
where
$$
C_{n,k}=\frac{(2k)!(n-3)!}{2^k(n-2k)!},
$$
$\nabla^*\nabla$ is the Bochner Laplacian acting on $\Cc^1(X)$, $\delta_\gm$ is the divergence operator, $\delta^*_\gm$ is the $L^2$ adjoint of $\delta_\gm$ and $\Delta_\gm=\delta_\gm d$ is the Laplacian acting on functions.

\begin{remark}\label{expdiv}
Our sign convention for $\delta_g$ is
$$
(\delta_g\omega)(x) = - \sum_i(\nabla_{e_i}\omega)(e_i,x),
$$
where $\nabla$ is the covariant derivative on $\Cc^1(X)$ induced by the Levi-Civita connection of $g$ so that $\Delta f=-f''$ on $\mathbb R$. Notice that we still denote by $\nabla$ the covariant derivative induced on the space of tensors of any kind.
\end{remark}

Thus, the identities $c_\gm\delta_{\gm}^*\delta_\gm h=-\delta_\gm\delta_\gm h$ and $c_\gm\nabla df=-\Delta_\gm f$ immediately give
\begin{eqnarray*}
c_\gm\dot{\mathcal R}_{\gm}^{(2k)}h & = & C_{n,k}\mu^{k-1}\bigg(n(k-1)\left(\Delta_\gm \mathsf {tr}_\gm h+\delta_\gm \delta_\gm h+ (n-1)\mu{\mathsf {tr}}_\gm h\right) +  \\
& &\qquad +({n-2k}) \Big(\Delta_\gm{\mathsf {tr}}_\gm h+\delta_\gm\delta_\gm h+(n-1)\mu {\mathsf {tr}}_\gm h\Big)-\\
& & \qquad \qquad -(n-2)(2k-1)(n-1)\mu{\mathsf {tr}}_\gm h\bigg)\\
& = & (n-2)C_{n,k}\mu^{k-1}\left(k(\Delta_\gm \mathsf {tr}_\gm h+\delta_\gm \delta_\gm h)-(n-1)(k-1)\mu{\mathsf {tr}}_\gm h)\right).
\end{eqnarray*}
On the other hand, $\mathcal R^{(2k)}_\gm=\lambda_k\gm$, where
$$
\lambda_k=\frac{(2k)!(n-1)!}{2^k(n-2k)!}
$$ by (\ref{lambda}). Thus, $(\dot c_\gm h){\mathcal R}_\gm ^{(2k)}=\lambda_k(\dot c_\gm h)\gm$. To compute this, we linearize the identity $c_gg=n$ to obtain $(\dot c_gh)g=-c_gh=-{\mathsf {tr}}_gh$, so that
$$
(\dot c_\gm h){\mathcal R}_\gm^{(2k)}=-C_{n,k}\mu^{k-1}(n-1)(n-2)\mu{\mathsf {tr}}_\gm h.
$$
The conclusion is that
\begin{equation}\label{div}
\dot{\mathcal S}_{\gm}^{(2k)}h=D_{n,k,\mu}\left(\Delta_\gm \mathsf {tr}_\gm h+\delta_\gm \delta_\gm h-(n-1)\mu{\mathsf {tr}}_\gm h\right),
\end{equation}
where
$$
D_{n,k,\mu}=\frac{1}{(2k)!}(n-2)kC_{n,k}\mu^{k-1}.
$$
Notice that for $k=1$ this gives
$$
2\dot{\mathcal S}_{\gm}^{(2)}h=\Delta_\gm \mathsf {tr}_\gm h+\delta_\gm \delta_\gm h-(n-1)\mu{\mathsf {tr}}_\gm h,
$$
a classical result \cite{B}.

In the sequel we shall use  a special case of (\ref{div}), which is obtained by setting $h=f\gm$, where $f\in\Omega^0(X)$. Since $\delta_\gm (f\gm)=-df$, we get
\begin{equation}\label{div2}
\dot{\mathcal S}_{\gm}^{(2k)}(f\gm)=D'_{n,k,\mu}\Lc_{\gm}f,
\end{equation}
where $D'_{n,k,\mu}=(n-1)D_{n,k,\mu}$ and
\begin{equation}\label{op}
\Lc_{\gm}=\Delta_{\gm}-n\mu.
\end{equation}

The following well-known result concerning the operator $\Lc_\gm$ above will play a key role in our analysis.

\begin{lemma}\label{specbotton}
If $(X,\gm)\in \Fc_n$ then the equation $\Lc_\gm w=0$ has no nontrivial solutions. In fact, $\lambda_1(\Delta_{\gm})$, the first positive eigenvalue of $\Delta_\gm$, is strictly greater than $n\mu$.
\end{lemma}

This is clearly the case if $\mu<0$ because $\Delta_\gm$ is non-negative, but if $\mu>0$ then the equation $\Lc_\gm w=0$ might have nontrivial solutions.
For instance, if $(X,\gm)$ is the round sphere with curvature $\mu$
then $\lambda_1(\Delta_\gm)=n\mu$  and the equation has a solution space of dimension $n+1$. But if $(X,\gm)$ is
a {\em nontrivial} quotient of a round sphere,  a theorem by Lichnerowicz and Obata \cite{BGM} implies   $\lambda_1(\Delta_\gm)>n\mu$, as desired.

\begin{remark}\label{labbi}
We note that a formula for $\dot\Sc^{(2k)}_g$, $g$ a general metric, has been obtained earlier \cite{Lo}, \cite{L2}. However, the main motivation of these authors is to extremize the  Hilbert-Einstein-Lovelock action (\ref{hillove}). This naturally led them to express  $\dot\Sc^{(2k)}_gh$ as a sum of a zero order term in $h$, which effectively contributes to the corresponding Euler-Lagrange equations, and a
rather cumbersome higher order term written in divergence form which vanishes after integration. Our computation above shows that, in the constant curvature case, the divergence term gives rise  to  the first two terms in the right-hand side of (\ref{div}).
\end{remark}

\section{The proof of Theorem \ref{main2}}\label{main2sec}

Let $(X,\gm)\in \Fc_n$ with sectional curvature $\mu\neq 0$. By scaling we may assume that $(X,\gm)\in\Mc_1(X)$, the space of unit volume metrics on $X$. Thus we must prove Theorem \ref{main2} with $\nu=1$.

In what follows the symbol $H_{\gm}^r(\Uc)$ denotes the standard Sobolev construction applied to a subset $\Uc$ of sections of a vector bundle over $X$, so that for instance $H_{\gm}^r(\Mc(X))$ is the Hilbert manifold, modeled on $H^r_\gm(\Cc^1(X))$, of metrics with derivatives up to order $r$ defined almost everywhere and square integrable (with respect to $\gm$).

We choose $r>\frac{n}{2}+4$ and define $\Ac_r:H^r_\gm(\Mc(X))\to H^{r-4}_\gm(\Omega^0_{\bullet}(X))$ by
$$
\Ac_r(g)=\Delta_{g}\Sc^{(2k)}_{g}-\int_X \Delta_{g}\Sc^{(2k)}_{g}\nu_\gm.
$$
Here,
$
\Omega^0_{\bullet}(X)=\left\{ \rho\in \Omega^0(X);\int_X \rho\,\nu_\gm=0\right\}.
$
Since $g\in H^r_\gm(\Mc(X))$ implies $R_{g}\in H^{r-2}_\gm(\Cc^2(X))$, $\Ac_r$ is well-defined and smooth due to the local expression (\ref{fully}) for $\Sc_{g}^{(2k)}$ and the fact that for $r-2>n/2+2>n/2$ the Sobolev space $H^{r-2}_\gm$ is a Banach algebra under pointwise multiplication \cite{MS}.

\begin{lemma}\label{lemma0}
There exists  a neighborhood, say $V^r$, of $\gm$ in $H_{\gm}^r(\Mc_1^k(X))$ which is a smooth submanifold of $H_{\gm}^{r}(\Mc(X))$ with $T_{\gm}V^r=\ker \dot\Ac_r(\gm)$.
\end{lemma}

For the proof first note that  $\dot\Delta_\gm (h)\Sc^{(2k)}_{\gm}=0$ for $h\in\Cc^1(X)$ because  $\Sc^{(2k)}_{\gm}$ is constant. Thus we get from (\ref{div}),
\begin{eqnarray}
\dot\Ac_r(\gm)(h) & = & \Delta_\gm\dot\Sc^{(2k)}_{\gm}(h)\nonumber\\
\label{nost} & = & D_{n,k,\mu}\Delta_\gm\left(\Delta_\gm \mathsf {tr}_\gm h+\delta_\gm \delta_\gm h-(n-1)\mu{\mathsf {tr}}_\gm h\right).
\end{eqnarray}
We now take $h=f\gm$ for $f\in H_{\gm}^{r}(\Omega^0(X))$, so that by (\ref{div2})
$$
\dot\Ac_r(\gm)(f\gm)=D'_{n,k,\mu}\Delta_\gm\Lc_\gm f.
$$
Since, by Lemma \ref{specbotton}, $n\mu$ is not an eigenvalue of $\Delta_\gm$, it follows that
$$
\dot\Ac_r(\gm):H_{\gm}^r(\Cc^1(X))\to H_{\gm}^{r-4}(\Omega^0_\bullet(X))
$$
is onto. The lemma is now an immediate consequence of the Implicit Function Theorem and the fact that $\Ac^{-1}_r(0)=\Mc_1^{k}(X)$.

\begin{lemma}\label{lemma}
If $\xi^r:H_{\gm}^r(\Omega^0_+(X))\times V^r\to H_{\gm}^r(\Mc(X))$ is the smooth map given by  $\xi^r(f,g)=fg$ then $d\xi^r_{(1,\gm)}$ is an isomorphism.
\end{lemma}

If $d\xi^r_{(1,\gm)}(\phi,h)=h+\phi \gm=0$ then $h=-\phi \gm\in \ker \dot\Ac_r(\gm)$ so that $\Lc_\gm\Delta_\gm\phi=0$. Hence $\Delta_\gm\phi=0$ by Lemma \ref{specbotton} and  $\phi$ is constant. But $\int_X\phi\,\nu_\gm=0$ because $V^r\subset H_{\gm}^r(\Mc_1(X))$ and hence $\phi=0$, which implies $h=0$. This shows the injectivity of $d\xi_{(1,\gm)}$.

As for the surjectivity, the identity
$$
{\rm Im}\,d\xi^r_{(1,\gm)}=T_\gm V^r\oplus H_{\gm}^r(\Omega^0(X))\gm
$$
already shows that ${\rm Im}\,d\xi^r_{(1,\gm)}$ is closed in $H_{\gm}^r(\Cc^1(X))$. Thus assume by contradiction the existence of $h\neq 0$ in $H_{\gm}^r(\Cc^1(X))$ orthogonal to both $T_\gm V^r$ and $H_{\gm}^r(\Omega^0(X)\gm$. As it is clear from (\ref{nost}), $\dot\Ac_r(\gm)$ has surjective symbol, and since $T_{\gm}V^r=\ker \dot\Ac_r(\gm)$ one has the decomposition \cite{E}
$$
H_{\gm}^r(\Cc^1(X))=\mathbb{R}\gm\oplus T_\gm V^r\oplus {\rm Im}\,\dot\Ac_r(\gm)^*,
$$
where $\dot\Ac_r(\gm)^*$ is the $L^2$ adjoint of $\dot\Ac_r(\gm)$. This allows us to write $h=\dot\Ac_r(\gm)^*(\varphi)$, that is,
$$
h=D_{n,k,\mu}\left((\Delta^2_\gm\varphi)\gm+\nabla d\Delta_\gm\varphi-(n-1)\mu (\Delta_\gm\varphi) g_\mu\right),
$$
so that by taking traces,
$$
{\mathsf {tr}}_\gm h=D'_{n,k,\mu}\Lc_\gm\Delta_\gm\varphi.
$$
Now, $\int_X{\mathsf {tr}}_\gm h\,\nu_\gm=0$ because $h$ is orthogonal to $H_{\gm}^r(\Omega^0(X))\gm$,
and using Lemma \ref{specbotton}, $\int_X\Delta_\gm\varphi\,\nu_\gm=0$ and the well-known variational characterization of  $\lambda_1(\Delta_\gm)$,
$$
n\mu <\lambda_1(\Delta_\gm)\leq \frac{\int_X|\nabla\Delta_{\gm}\varphi|^2\nu_{\gm}}
{\int_X |\Delta_{\gm}\varphi|^2\nu_{\gm}}= n\mu,
$$
a contradiction unless
$\Delta_\gm\varphi=0$ i.e. $\varphi$ is constant and hence $h=0$. This completes the proof of Lemma \ref{lemma}.

With Lemmas \ref{lemma0} and \ref{lemma} at hand, it is now straightforward to carry out the proof of Theorem \ref{main2}, by essentially using the fact that objects in the ILH category are defined as inverse limits of objects in the $H_{\gm}^r$ category as $r\to +\infty$.  We will omit the details and instead will refer to \cite{K}; see the proof of Theorem 2.5 there.

\begin{remark}\label{sum}
We note that the arguments above can be easily adapted to the  Yamabe  problem for  certain functions of the Gauss-Bonnet curvatures. Indeed, if $(X,\gm)\in\Fc_n$ let $k_n$ be the largest integer strictly less than $n/2$. For $G:\mathbb R^{k_n}\to \mathbb R$  smooth define
\begin{equation}\label{linearcomb}
\Gc_g=G(\Sc^{(2)}_{g},\ldots,\Sc^{(2k_n)}_{g}),\qquad g\in\Mc(X),
\end{equation}
and {\em assume} that, throughout $(X,\gm)$,
\begin{equation}\label{sumassump}
D\equiv\sum_{1\leq k\leq k_n }D_{n,k,\mu}\frac{\partial G}{\partial x_k}\left(\Sc^{(2)}_{g_\mu},\ldots,\Sc^{(2k_n)}_{g_\mu}\right)\neq 0.
\end{equation}
In this way we get the formulae corresponding to (\ref{div}) and (\ref{div2}), namely,
$$
\dot\Gc_\gm(h)=D\left(\Delta_\gm \mathsf {tr}_\gm h+\delta_\gm \delta_\gm h-(n-1)\mu{\mathsf {tr}}_\gm h\right),
$$
and
$$
\dot\Gc_\gm(fg)=D'\Lc_\gm f,\qquad D'=(n-1)D,
$$
so that a result similar to Theorem \ref{main2} holds: there exists a neighborhood $U$ of $\gm$ in $\Mc(X)$ such that any metric in $U$ is conformal to a metric $g$ with  $\Gc_{g}$ a constant.
\end{remark}

\end{document}